\documentclass[a4paper,12pt]{amsart}
\usepackage{amsfonts,amsthm,amssymb}
\usepackage{epsfig}
\usepackage{pstricks,slashbox,multirow}
\usepackage{graphicx}
\usepackage{amsmath}
\usepackage{latexsym}
\usepackage{tikz}
\usepackage{fancyvrb}
\usepackage{datetime}
\usepackage{verbatim}
\newtheorem{theorem}{Theorem}[section]

\newtheorem{definition}[theorem]{Definition}

\newtheorem{exm}{Example}

\newtheorem{res}[theorem]{Result}
\newtheorem{note}[theorem]{Note}

\title[A study on edge coloring and edge sum coloring of integral sum graphs]{A study on edge coloring and edge sum coloring of integral sum graphs}
\author{\sc V. Vilfred Kamalappan} 
\address{Department of Mathematics, ~Central ~University ~of~ Kerala,~Periye, \linebreak  ~Kasaragod,~ Kerala, ~India~ - ~671 316.}
\email{vilfredkamalv@cukerala.ac.in}
\author{\sc Lowell W. Beineke} 
\address{Purdue~ University, Fort~ Wayne, Indiana~ 46805, U.S.A.}
\email{beineke@pfw.edu}
\author{\sc L. Mary Florida} 
\address{  St. Xavier’s ~Catholic ~College~ of~ Engineering, ~Chunkankadai,  \linebreak Nagercoil, Kanyakumari District, Tamil ~Nadu,  India - 629 807.}
\email{link2florida@yahoo.co.in}
\author{\sc Julia K. Abraham} 
\address{Department of Mathematics, ~Central ~University ~of~ Kerala,~Periye, \linebreak  ~Kasaragod,~ Kerala, ~India~ - ~671 316.}
\email{juliakabraham92@gmail.com}

\thanks{Dedicated to the memory of Professor Frank Harary  on his Centenary Year}

\subjclass[2010]{05C78, 05C15, 05C75.}
\keywords{Integral sum graph, Star graph, edge sum number, edge-sum class, edge sum color graph, edge sum chromatic number.}

\date{}

\begin{document}

\begin{abstract} Frank Harary introduced the concept of integral sum graph.  A graph $G$ is an \emph{ integral sum graph}  if its vertices can be labeled with distinct integers so that $e = uv$ is an edge of $G$ if and only if the sum of the labels on vertices $u$ and $v$ is also a label in $G.$ For any non-empty set of integers $S$, let $G^+(S)$ denote the integral sum graph on the set $S$. In $G^+(S)$, we define an \emph{edge-sum class} as the set of all edges each with same edge sum number and call $G^+(S)$ an \emph{edge sum color graph} if each edge-sum class is considered as an edge color class of $G^+(S)$. The number of distinct edge-sum classes of $G^+(S)$ is called its \emph{ edge sum chromatic number}. The main results of this paper are (i) the set of all edge-sum classes of an integral sum graph partitions its edge set; (ii) the edge chromatic number and the edge sum chromatic number are equal for the integral sum graphs $G_{0,s}$ and $S_n$, Star graph of order $n$, whereas it is not in the case of $G_{r,s} = G^+([r,s])$, $r < 0 < s$, $n,s \geq 2$, $n,r,s\in\mathbb{N}$. We also obtain an interesting integral sum labeling of Star graphs.
\end{abstract}
	
\maketitle

\section{Introduction}
	
In 1994 Frank Harary \cite{ha90}, \cite{ha94} introduced the concepts of sum and integral sum graphs. A graph $G$ is a \emph{ sum graph} or \emph{  $\mathbb{N}$$-$$sum$ graph} if its vertices can be labeled with distinct positive integers so that $e = uv$ is an edge of $G$ if and only if the sum of the labels on vertices $u$ and $v$ is also a label in $G.$ An \emph{ integral sum graph} or \emph{ $\mathbb{Z}$-sum graph} is defined similarly, the only difference being that the labels can be integers. For any non-empty set of integers $S$, let $G^+(S)$ denote the integral sum graph on the set $S$. Frank Harary also produced a family of integral sum graphs $G_{-n, n}$ = $G^+([-n, n])$ which is generalized to $G_{r, s}$ where $[r, s]$ = $\{r,r + 1,\ldots,s\},$ $r,s,n\in \mathbb{Z}$ in \cite{vmf12}.  Various extensions to the notion of sum graphs were introduced by several authors \cite{ga20, lhv}. 

The study of sum graphs is used to distribute secret information to a set of people so that only an authorized set of people can reconstruct the secret \cite{ssm} and also used to storage and manipulation of relational database \cite{ga20}.

	We adopt some notation that we hope will assist in keeping formulas relatively brief:
	
	\begin{enumerate}
	  \item [\rm 1.] $n$ will always denote a positive integer, and $G_n$ is the sum graph $G^+([1,n])$.
		
		\item [\rm 2.] $r$ will always denote a negative integer and $s$ a positive integer with $r + s \geq 0$, and $G_{r,s}$ denotes the integral sum graph $G^+([r,s])$. 
		\item [\rm 3.] The number of vertices in a graph $G$ will be denoted by $|G|$  and the number of edges by $||G||$.
	\end{enumerate}
	
	We note that $G^+([-s,-1]) \cong G^+([1,s])$ so that a labeling with only negative labels is the same as one with only positive labels.  More generally, if every label in a sum graph is replaced by its negative, then the two graphs are isomorphic, so the condition $r + s \geq 0$ on $G_{r,s}$ is not really a restriction at all.
	
	Properties of sum and integral sum graphs have been investigated by various authors  \cite{ch90}-\cite{ga20},\cite{ha90}-\cite{wl15}. It is noted in \cite{vmf12} that the interval graphs $G_{r,s}$ can be expressed in terms of the graphs $G_n$ as $G_{r,s}$ $\cong$ $K_1 * G_{-r} * G_s$. Note that the graph operation of the join, which we denote here by $*$, is both associative and commutative. Integral sum graphs $G_{0,6}, G_{-1,5}, G_{-2,4},$ $G_{-3,3}$ are given in Figures 1-4.

	\vspace{.5cm}	
	\begin{center}
	\begin{tikzpicture}
	
	\node (a6) at (9,1.75)  [circle,draw,scale=0.8] {6};
	\node (a0) at (11,2.75)  [circle,draw,scale=0.8]{0};
	\node (a1) at (13,1.75)  [circle,draw,scale=0.8]{1};
	\node (a2) at (13.5,.5)  [circle,draw,scale=0.8] {2};
	\node (a3) at (12,-.5)  [circle,draw,scale=0.8] {3};
	\node (a4) at (10,-.5)  [circle,draw,scale=0.8] {4};
	\node (a5) at (9,.5)  [circle,draw,scale=0.8]{5};

	\draw (a0) -- (a1);
	\draw (a0) -- (a2);
	\draw (a0) -- (a3);
	\draw (a0) -- (a4);
	\draw (a0) -- (a5);
	\draw (a0) -- (a6);
	
	\draw (a1) -- (a2);
	\draw (a1) -- (a3);
	\draw (a1) -- (a4);
	\draw (a1) -- (a5);
	
	\draw (a2) -- (a3);
	\draw (a2) -- (a4);
	\node (b6) at (15,1.75)  [circle,draw,scale=0.8] {-1};
	\node (b0) at (17,2.75)  [circle,draw,scale=0.8]{0};
	\node (b1) at (19,1.75)  [circle,draw,scale=0.8]{1};
	\node (b2) at (19.5,.5)  [circle,draw,scale=0.8] {2};
	\node (b3) at (18,-.5)  [circle,draw,scale=0.8] {3};
	\node (b4) at (16,-.5)  [circle,draw,scale=0.8] {4};
	\node (b5) at (15,.5)  [circle,draw,scale=0.8]{5};

	\draw (b0) -- (b1);
	\draw (b0) -- (b2);
	\draw (b0) -- (b3);
	\draw (b0) -- (b4);
	\draw (b0) -- (b5);
	\draw (b0) -- (b6);
	
	\draw (b6) -- (b1);
	\draw (b6) -- (b2);
	\draw (b6) -- (b3);
	\draw (b6) -- (b4);
	\draw (b6) -- (b5);
		
	\draw (b1) -- (b2);
	\draw (b1) -- (b3);
	\draw (b1) -- (b4);
		
	\draw (b2) -- (b3);
	
	\end{tikzpicture}

\vspace{.2cm}		
Fig. 1. $G_{0, 6}$ \hspace{4cm} Fig. 2. $G_{-1, 5}$
\end{center}
	
	\vspace{.4cm}	
	\begin{center}
	\begin{tikzpicture}
	
	\node (c6) at (17,2.25)  [circle,draw,scale=0.8] {-1};
	\node (c0) at (19,3)  [circle,draw,scale=0.8]{0};
	\node (c1) at (21,2.25)  [circle,draw,scale=0.8]{1};
	\node (c2) at (22,.75)  [circle,draw,scale=0.8] {2};
	\node (c3) at (21,-.75)  [circle,draw,scale=0.8] {3};
	\node (c4) at (19,-.75)  [circle,draw,scale=0.8] {4};
	\node (c5) at (17,.75)  [circle,draw,scale=0.8]{-2};
	
	\draw (c0) -- (c1);
	\draw (c0) -- (c2);
	\draw (c0) -- (c3);
	\draw (c0) -- (c4);
	\draw (c0) -- (c5);
	\draw (c0) -- (c6);
	
	\draw (c6) -- (c1);
	\draw (c6) -- (c2);
	\draw (c6) -- (c3);
	\draw (c6) -- (c4);
	
	\draw (c5) -- (c1);
	\draw (c5) -- (c2);
	\draw (c5) -- (c3);	
	\draw (c5) -- (c4);
	
	\draw (c1) -- (c2);
	\draw (c1) -- (c3);	
		
	\node (d6) at (24,2.25)  [circle,draw,scale=0.8] {-1};
	\node (d0) at (26,3)  [circle,draw,scale=0.8]{0};
	\node (d1) at (28,2.25)  [circle,draw,scale=0.8]{1};
	\node (d2) at (28,.75)  [circle,draw,scale=0.8] {2};
	\node (d3) at (27.5,-.75)  [circle,draw,scale=0.8] {3};
	\node (d4) at (24.5,-.75)  [circle,draw,scale=0.8] {-3};
	\node (d5) at (24,.75)  [circle,draw,scale=0.8]{-2};

	\draw (d0) -- (d1);
	\draw (d0) -- (d2);
	\draw (d0) -- (d3);
	\draw (d0) -- (d4);
	\draw (d0) -- (d5);
	\draw (d0) -- (d6);
	
	\draw (d6) -- (d1);
	\draw (d6) -- (d2);
	\draw (d6) -- (d3);
	\draw (d6) -- (d5);
	
	\draw (d5) -- (d1);
	\draw (d5) -- (d2);
	\draw (d5) -- (d3);
		
	\draw (d4) -- (d1);
	\draw (d4) -- (d2);
	\draw (d4) -- (d3);
	
	\draw (d1) -- (d2);
			
	\end{tikzpicture}
	
	\vspace{.2cm}
	Fig. 3. $G_{-2, 4}$ \hspace{5cm} Fig. 4. $G_{-3, 3}$
\end{center}
	\vspace{.5cm}

	The following results will be useful in subsequent sections.
	
	\begin{theorem} {\rm{\cite{vmf12}} \label{a1} 
	 For integers $r$ and $s$ with $r < 0 < s$, 	$G_{r,s}$ $\cong$ $K_1 * G_{-r} * G_s$. \hfill $\Box$}
	\end{theorem}	
	
	\begin{theorem}  {\rm{\cite{vmf12}} \label{a2}
		The degree of the vertex with label $i$ in $G_n$ is 
		\[ \deg (i) = \left\{ \begin{array}{ll}
		n-i-1 & \mbox{if $1 \leq i \leq \lfloor\frac{n}{2}\rfloor$} \\
		n-i & \mbox{if $\lfloor\frac{n}{2}\rfloor < i \leq n$. ~  }
		\end{array}  \right.  \] ~ \hfill $\Box$}
	\end{theorem}	
	Thus, the degree sequence of $G_{2k}$ is $\{2k-2, 2k-3, \ldots, k-1,$ $k-1,$ $k-2, ..., 1, 0\}$, while that of $G_{2k+1}$ is $\{2k-1, 2k-2, ..., k, k, k-1, ..., 1, 0\}$. By adding these numbers, we find the number of edges in $G_n$.  The graph $G_{2k}$ has $k(k-1)$ edges and the graph $G_{2k+1}$ has $k^2$ edges.  We combine these results in the next theorem.
	\begin{theorem}  {\rm{\cite{vbs}} \label{a3} 
		The graph $G_n$	has $\left\lfloor\frac{(n-1)^2}{4}\right\rfloor$ edges, $n\in\mathbb{N}$. ~ \hfill $\Box$}
	\end{theorem}	
	
	We now turn to the degrees of the vertices in the integral sum graphs $G_{r,s}$ (with our convention that $r < 0 < s$).
	
	\begin{theorem}  {\rm{\cite{vbs}} \label{a4}
		The degree of the vertex with label $i$ in $G_{r,s}$, with $n =  s - r + 1$, is 
		
		\[ \deg (i) = \left\{ \begin{array}{ll}
			n-i-1 & \mbox{if $1 \leq i \leq \lfloor\frac{~s~}{2}\rfloor$} \\
			n-i & \mbox{if $\lfloor\frac{~s~}{2}\rfloor < i \leq s$ }\\
		n-1  & \mbox{if $i = 0$}\\
		n+i-1  & \mbox{if $1 \leq -i \leq \left\lfloor \frac{-r}{2}\right\rfloor $} \\
		n+i & \mbox{if $ \left\lfloor \frac{-r}{2}\right\rfloor  < -i \leq -r $.}\\
\end{array} 	 \right.
		\] ~ \hfill $\Box$}
	\end{theorem}	
	
	By combining Theorems \ref{a1} and \ref{a3}, we find $m(r, s) = ||G_{r,s}||$, the number of edges in $G_{r,s}$ (recall that $r$ is negative) as follows.
	
\begin{theorem}  {\rm{\cite{vbs}} \label{a5}
		The number of edges in 	$G_{r,s}$ is 
		
$||G_{r,s}|| = -rs -r + s + \left\lfloor\frac{(r+1)^2}{4}\right\rfloor + \left\lfloor\frac{(s-1)^2}{4}\right\rfloor.$ 

\hspace{1.13cm}	$ = \frac{1}{4}(r^2+s^2-3r +3s-4rs)- \frac{1}{2}(\left\lfloor \frac{-r}{2}\right\rfloor + \left\lfloor \frac{s}{2} \right\rfloor)$.
\hfill $\Box$ }
\end{theorem}	
	
	\vspace{.2cm}
	There may be advantages to expressing this in odd and even cases:
	\[||G_{r,s}|| = \left\{ \begin{array}{ll}
	a^2 + b^2 + 4ab + 4a + 4b + 3 \\ \hspace{4cm} \mbox{if $r = -(2a+1)$ and $s = 2b+1$}\\ 
	a^2 + b^2 + 4ab + 2a + 3b + 1 \hspace{.1cm} \mbox{if $r = -(2a+1)$ and $s = 2b$}\\ 
	a^2 + b^2 + 4ab + 3a + 2b + 1 \hspace{.1cm} \mbox{if $r = -2a$ and $s = 2b+1$}\\ 
	a^2 + b^2 + 4ab + a + b \hspace{1.25cm} \mbox{if $r = -2a$ and $s = 2b$.}\\ 
	\end{array}  \right.
	\]		
	
	In considering different sum graphs and their subgraphs, we generally disregard their labels, as in the following definitions.  In general, two graphs are said to be \emph{comparable} if either is a subgraph of the other, and \emph{non-comparable} otherwise. 
	
	\emph{Star graph} $S_{n}$ on $n$ vertices is defined as $S_{n}$ = $K_1 * ((n-1) K_1)$ = $K_{1, n-1}$, $n \geq 2$ and $n\in\mathbb{N}$. Star graphs are integral sum graphs \cite{ch90}. The following result presents a new and more general integral sum labeling of Star graphs and the sequence of labeling an interesting one. 

\begin{res}\quad {\rm \label{a51} Let $S_n$ be a star graph, $V(S_n)$ = $\{u_0, v_1, v_2, . . . ,$ $v_{n-1}\}$, $d(u_0) = n-1$ and $d(v_1)$ = 1 = $d(v_2)$ = . . . = $d(v_{n-1})$, $n \geq 2$ and $d,i,j,n,t\in\mathbb{N}$. Define labeling $f: V(S_n) \to \mathbb{Z}$ $\ni$ 

 $f(u_0) = 0$, $f(v_1) = t$, 

$f(v_{i+1})$ = $d\left(\sum^i_{j=1}f(v_j)\right) + t$, $i$ = 1,2,...,$n-1$. 

Show that $f$ is an integral sum labeling of $S_n$ and the sequence of vertex labeling of $v_i$s is $\{t{(d+1)}^{i-1}\}^n_{i=1}$. }
\end{res}
 
\begin{proof}\quad Using induction on $i$, it is easy to prove that $f(v_{i+1})$ = $t{(d+1)}^{i}$, $i$ = 1,2,...,$n-1$ and $d,t\in\mathbb{N}$. 

Also, (i)~ for $i = 1,2,...,n$, $f(u_0) + f(v_i) = 0 + t{(d+1)}^{i-1}$ = $f(v_i)$ and (ii)~ for $1 \leq i < j < k \leq n$ and $i,j,k,d,t\in\mathbb{N}$, $f(v_i) + f(v_j)$ $\neq$ $f(v_k)$ since $t{(d+1)}^{i}$ + $t{(d+1)}^{j}$ = $t{(d+1)}^{i}(1+ {(d+1)}^{j-i})$ $\neq$ $t{(d+1)}^{k}$, $1 \leq j-i < j < k$. This implies that $u_0$ and $v_i$ are adjacent for every $i$ whereas $v_j$ and $v_k$ are non-adjacent for every $j$ and $k$, $j \neq k$, $1 \leq j,k \leq n$. Thus, $f$ is an integral sum labeling of $S_n$, $n \geq 2$.
\end{proof}

\begin{note}\quad {\rm The above sequence $\{t{(d+1)}^{i-1}\}^n_{i=1}$ of vertex labeling $f(v_i)$s of $S_n$ is a GP and also is a strictly monotonic increasing sequence for a given value of $t$ and $d$, $n \geq 2$ and $d,i,n,t\in\mathbb{N}$.}
\end{note}

\begin{definition}\quad	A sum graph is \emph{sum-maximal} if it is not isomorphic to a proper spanning subgraph of another sum graph and is a \emph{sum-maximum} graph if it has the greatest number of edges among all sum-maximal graphs of a given order.
\end{definition}
 For example, we observe that for $s \geq 2$, $G_{0,s}$ is a proper subgraph of $G_{-1,s-1}$, so that $G_{0,s}$ is not sum-maximal for $s \geq 2$ whereas $G_{r,s}$ is sum-maximal of order $s-r+1$ for $r < 0 < s$ and $ s \geq - r$ \cite{vmf12}. The results in the following two theorems are related to the above properties on $G_{r,s}$.
	
	\begin{theorem} {\rm{\cite{vmf12}} \label{a6}
		For $r < 0 < s$ and $ s \geq - r$, 	$G_{r,s}$ is sum-maximal of order $s-r+1$. \hfill $\Box$}
	\end{theorem}	
	
	\begin{theorem} {\rm{\cite{vmf12}} \label{a7}
		For  $n \geq 2$, $G_{- \lfloor \frac{n}{2} \rfloor, \lceil \frac{n}{2} \rceil}$ is a sum-maximum graph of order $n+1$.  \hfill  $\Box$}
	\end{theorem}	

In this paper, we define edge-sum class, edge sum color graph and edge-sum chromatic number of integral sum graphs and prove that (i) the set of all edge-sum classes of an integral sum graph partitions its edge set; (ii) the edge chromatic number and the edge sum chromatic number are equal for the integral sum graphs $G_{0,s}$ and $S_n$ whereas it is not in the case of $G_{r,s} = G^+([r,s])$, $r < 0 < s$, $n,s \geq 2$, $n,r,s\in\mathbb{N}$. 

Studying edge coloring of integral sum graphs $G_{r,s}$ is the motivation to do this work. 
	
\section{On edge sum coloring of integral sum graphs} 

	We begin with some notation.  Given an integral sum graph $G^+(S)$, we assume that vertex $u_j$ has label $j$. Every edge of an integral sum graph has an induced \emph{edge sum number} which is the sum of the labels of its end vertices and adjacent edges have different edge sum numbers. In an integral sum graph, the set of all edges each with same edge sum number, say $i$, is called the \emph{edge-sum class} and is denoted by $E_i$. For example, consider $G^+(S) = G_{-2, 4}$ as shown in Fig. 3. Then we have these edge-sum classes:
	
	$E_{-2} = \{(-2,0)\}$,
	
	$E_{-1} =\{(-1,0),(-2,1)\}$,
	
	$E_{0}=\{(-2,2), (-1,1)\}$,
	
	$E_{1} = \{(-2, 3), (-1,2), (0, 1)\}$,  
	
	$E_{2} = \{(-2, 4), (-1,3), (0, 2)\}$, 
	
	$E_{3} = \{(-1,4), (0,3), (1, 2)\}$ and
	
	$E_{4} = \{(0,4), (1, 3)\}$.
	
	The reason for calling $E_i$ as the edge-sum class is based on Theorems \ref{c1} and \ref{c2}.
	
	\vspace{.2cm}
We now consider the sum graphs what one vertex $i$ of a given sum graph contributes, and that is not just the vertex $i$ itself and its incident edges, but also all edges with sum $i$.  It is therefore convenient to define the set $F_i = \{i\} \cup E_i$. 

		 In this section, we study interesting properties of edge-sum classes of integral sum graph $G^+(S)$, relates these with edge coloring and a few open problems at the end.  Recall that $E_k$ denotes the set of edges, each with edge sum $k$.
		
		\begin{theorem} {\rm \label{c1}
			Let $G^+(S)$  be an integral sum graph. Then the following holds in $G^+(S)$: 
			\begin{enumerate} 
				\item [\rm (a)]  $E(G^+(S)) = \bigcup\limits_{i \in S} E_i$. 
				\item [\rm (b)] Two edge-sum classes are either equal or disjoint.   
				\item [\rm (c)] For $i \neq j$ and $i,j \in S$ if $E_i,E_j \neq \emptyset$, then $E_i \neq E_j$. 
				\item [\rm (d)] A non-empty edge-sum class is an independent set of edges in the integral sum graph. 
				\item [\rm (e)] The number of distinct non-empty edge-sum classes of $G^+(S)$ is at the most the order of the graph, with equality if and only if every vertex label is the induced sum of some edge(s) in $G^+(S)$.
		\end{enumerate}}
		\end{theorem}
		\begin{proof} 
		\begin{enumerate} 
				\item [\rm (a)]  Every edge of $G^+(S)$ has an induced edge sum number, say $i$, which is the sum of the labels of its end vertices and $E_i$ is the set of all edges, each with induced edge sum number $i$  in $G^+(S)$, $i\in S$. Also, either $E_i$ = $\emptyset$, if there is no edge in $G^+(S)$ with induced edge sum number $i$ or $E_i$ $\neq$ $\emptyset$, if $G^+(S)$ has at least one edge whose edge sum number is $i$, $i\in S$. 
				
				Hence, $E(G^+(S)) = \bigcup\limits_{i \in S} E_i$.
			\item [\rm (b)]   The result is true when $E_i$ = $\emptyset$ or $E_j$ = $\emptyset$ or $E_i$ = $E_j$ = $\emptyset$, $i,j\in S$. 
			
			Let $E_i,E_j \neq \emptyset$, $i,j\in S$. Then, $i$ = $j$ if and only if $E_i$ = $E_j$ since $E_i$ is the set of all edges of $G^+(S)$ with induced edge sum number $i$, $i,j\in S$. Hence, we get $(b)$.
			\item [\rm (c)]   Follows from $(b)$.
			\item [\rm (d)]  Let $e_i = uv$ and $e_j = uw$ be two adjacent edges in $G^+(S)$, $i \neq j$ and $i,j\in S$. By the definition of integral sum labeling, vertices $v$ and $w$ can not have same integral sum labeling in $G^+(S)$ and thereby induced edge sum numbers of $uv$ and $uw$ are different. And so $e_i$ and $e_j$ belong to different edge-sum classes. Hence, we get $(d)$. 
			\item [\rm (e)]  Follows from the fact that every $i\in S$ need not be an induced edge sum number of an edge in $G^+(S)$ but every induced edge sum number belongs to $S$.
\end{enumerate}
\end{proof}
		
		\begin{theorem} {\rm \label{c2}
			Let $G^+(S)$ be an integral sum graph with $E(G^+(S)) \neq \emptyset$. Then the set of all edge-sum classes of $G^+(S)$ partitions its edge set.}
		\end{theorem}
		\begin{proof} \quad In an integral sum graph $G^+(S)$, every edge has an induced edge sum number which belongs to $S$ and also,  $E(G^+(S)) = \bigcup\limits_{i \in S} E_i$ Theorem \ref{c1}.(a). By $(b)$ and $(c)$ of Theorem \ref{c1}, any two non-empty edge-sum classes in an integral sum graph are either equal or disjoint. Thus, we get the result.  
		\end{proof}
		
		We have seen that the set of all edge-sum classes of an integral sum graph partitions its edge set. For a given integral sum graph $G^+(S)$, the set of all edge-sum classes is unique. Property $(d)$ helps us to consider an integral sum graph as an \emph{ edge sum color graph} by applying same color to all edges in an edge-sum class and different colors to (edges of) different edge-sum classes. In an integral sum graph $G$ with vertex labeling function $f$ the induced edge sum number of an edge is the sum of the labels of its end vertices. From these, we define an edge sum color graph as follows. 
		\begin{definition} \label{c3}
			An integral sum graph $G$ is called an \it{edge sum color graph} if there exists a function $c:$ $V(G)~\rightarrow~ \mathbb{Z}$ such that any two edges of $G$ have the same color if and only if their edge sum numbers are same. 
\end{definition}			
			Clearly, in an edge sum color graph with respect to an integral sum labeling $f$, edges of an edge-sum class take same color and edges of different edge-sum classes have different colors. And edge sum color graph is an integral sum graph only. Hence, whenever we say graph $G$ is an edge sum color graph, then it is understood that the graph $G$ we consider is an integral sum graph.    
		
\begin{definition} \quad	In an integral sum graph $G$, the number of distinct (non-empty) edge-sum classes of $G$ is called the \emph{ edge sum chromatic number} of $G$ and is denoted by $\chi^{'}_{\mathbb{Z}-sum}(G)$.
\end{definition}
		
		It is clear that for a given integral sum graph $G^+(S),$ its edge sum chromatic number $\chi^{'}_{\mathbb{Z}-sum}(G^+(S))$ = order of the graph $G^+(S)$ if and only if $E_i \neq \emptyset$ for every $i \in S$.
		
		For an integral sum graph, its edge-sum classes are unique whereas its edge color classes need not be.  
		
		Also, in an integral sum graph $G$, if we consider each of its edge-sum class as an edge color class, then this edge coloring of $G$ may not be a minimal edge coloring of $G$. The following example illustrates the above.

\begin{exm} \quad For the integral sum graph $G_{-1,5}$, $\chi^{'}(G_{-1,5})$ = 6 $\neq$ $\chi^{'}_{\mathbb{Z}-sum}(G_{-1,5})$ = 7.
\end{exm}

 The edge-sum classes of $G_{-1,5}$ are 

$E_{-1}$ = $\{(-1,0)\},$ 

$E_0 = \{(-1,1)\},$ 

$E_{1} = \{(0,1),(-1,2)\},$ 

$E_{2} = \{(0,2), (-1,3)\},$ 

$E_{3} = \{(0,3), (1, 2), (-1,4)\},$ 

$E_{4}$ = $\{(0,4), (1,3), (-1,5)\}$ and 

$E_{5} = \{(0,5),$ $(1,4),(2,3)\}$. This implies, the edge sum chromatic number $\chi^{'}_{\mathbb{Z}-sum}(G_{-1,5})$ = 7. 

One set of edge color classes of $G_{-1,5}$ is given by 

$\{\{(0, 1), (-1, 2)\},$ (say gray color); 

$\{(0, 2), (-1, 3)\},$ (say orange color);

$\{(0, 3), (1, 2), (-1,4)\},$ (say cyan color);

$\{(0, 4), (1, 3), (-1,5)\},$ (say blue color);

$\{(0, 5), (-1, 1)\}$ and (say green color);

$\{(0, -1), (1, 4), (2, 3)\}\}$ (say brown color). This implies, the edge chromatic number $\chi^{'}(G_{-1,5})$ = 6, follows from Vizing's theorem \cite{ha69} since $\Delta (G_{-1,5}) = 6 = deg(0)$ = $deg(-1)$. 

Thus $G_{-1,5}$ has 7 different edge-sum classes whereas it’s edge chromatic number is 6. This implies that $\chi^{'}(G_{-1,5})$ $\neq$ $\chi^{'}_{\mathbb{Z}-sum}(G_{-1,5})$. 

In Fig. 5, edges of $G_{-1,5}$ are colored as given in the above example. See Figure 5.
 
\vspace{.5cm}	
	\begin{center}
	\begin{tikzpicture}
		
	\node (b6) at (15,1.75)  [circle,draw,scale=0.8] {-1};
	\node (b0) at (17,2.75)  [circle,draw,scale=0.8]{0};
	\node (b1) at (19,1.75)  [circle,draw,scale=0.8]{1};
	\node (b2) at (19.5,.5)  [circle,draw,scale=0.8] {2};
	\node (b3) at (18,-.5)  [circle,draw,scale=0.8] {3};
	\node (b4) at (16,-.5)  [circle,draw,scale=0.8] {4};
	\node (b5) at (15,.5)  [circle,draw,scale=0.8]{5};

	\draw (b0) [gray, thick] -- (b1);
	\draw (b0) [orange, thick] -- (b2);
	\draw (b0) [cyan, thick] -- (b3);
	\draw (b0) [blue, thick] -- (b4);
	\draw (b0) [green, thick]-- (b5);
	\draw (b0) [brown, thick] -- (b6);
	
	\draw (b6) [green, thick] -- (b1);
	\draw (b6) [gray, thick] -- (b2);
	\draw (b6) [orange, thick] -- (b3);
	\draw (b6) [cyan, thick] -- (b4);
	\draw (b6) [blue, thick] -- (b5);

	\draw (b1) [cyan, thick] -- (b2);
	\draw (b1) [blue, thick]  -- (b3);
	\draw (b1) [brown, thick] -- (b4);
		
	\draw (b2) [brown, thick] -- (b3);
	
	\end{tikzpicture}

\vspace{.2cm}		
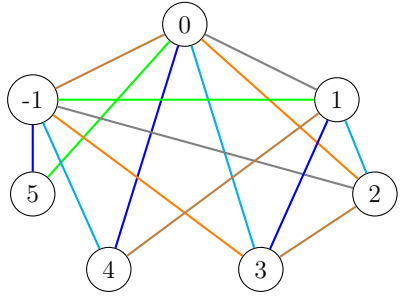
 Fig. 5. Edge colored $G_{-1, 5}$
\end{center}

	The following theorem presents relation between edge chromatic number and edge sum chromatic number of integral sum graphs (i)~ Star graphs $S_{n}$ on $n$ vertices, (ii)~ $G_{-1,1}$, (iii)~ $G_{0,s}$, $s\in \mathbb{N}$, (iv)~ $G_{-1,s}$ and $G_{-s,s}$, $s \geq 2$. 
		
		\begin{theorem} {\rm \label{c5}
For all negative integers $r$ and positive integers $s$, the following hold:
\par
\noindent{\rm {(a)}}   $\chi^{'}_{\mathbb{Z}-sum}(G_{r,s}) = |r|+s+1$;
\par
\noindent{\rm {(b)}}  $\chi^{'}(G_{-1,1}) = \chi^{'}_{\mathbb{Z}-sum}(G_{-1,1})$ = 3 and 

\hspace{.35cm} $\chi^{'}(G_{0,s}) = \chi^{'}_{\mathbb{Z}-sum}(G_{0,s})$ = $s$;

\par
\noindent{\rm {(c)}}  $\chi^{'}(S_{n}) = \chi^{'}_{\mathbb{Z}-sum}(S_{n})$ = $n-1$, $n \geq 2$; 
\par
\noindent{\rm {(d)}} $\chi^{'}(G_{-1,s}) \neq \chi^{'}_{\mathbb{Z}-sum}(G_{-1,s})$ for $s \geq 2$ and 
\par
\noindent{\rm {(e)}}  $\chi^{'}(G_{-s,s}) \neq \chi^{'}_{\mathbb{Z}-sum}(G_{-s,s})$ for $2 \leq s \leq 6$.}
		\end{theorem}

		\begin{proof}
		\begin{enumerate} 
\item [\rm {(a)}]    The edge-sum classes of $G_{r,s}$ are $E_r, E_{r+1}, \ldots, E_{-1}, E_0,$ $E_1,\ldots,E_s$ and each one is non-empty. Hence, $\chi^{'}_{\mathbb{Z}-sum}(G_{r,s})$ = $|r|+s+1$ for $-r,s \in \mathbb{N}$.\\
			\par
\noindent
\item [\rm {(b)}]  By observation, it is clear that  $\chi^{'}_{\mathbb{Z}-sum}(G_{-1,1}) = 3 =  \chi^{'}(G_{-1,1})$. 

For $s \in \mathbb{N}$, the edge-sum classes of $G_{0,s}$ are $E_0, E_1,\ldots,$ $E_s$ and each one, except $E_0$, is non-empty and hence $\chi^{'}_{\mathbb{Z}-sum}(G_{0,s})$ $= s$. Using Vizing's theorem \cite{ha69}, we get $s \leq$ $\chi^{'}(G_{0, s}) \leq s+1$ since $\Delta (G_{0, s}) = s = deg(0)$. Now, $\{\{(0,1)\},$ $\{(0,2)\},$ $\{(0,3), (1,2)\},$ $\{(0,4), (1,3)\},$ $\{(0,5), (1, 4),$ $(2, 3)\}, \{(0,6), (1, 5),$ $(2,4)\},$ $\ldots,$ $\{(0, s), (1,s-1),$ $(2,s-2), \ldots,$ $(\left\lfloor \frac{s-1}{2}\right\rfloor, \left\lfloor \frac{s+2}{2}\right\rfloor)\}\}$ is a set of edge color classes of $G_{0,s}$ and it is of order $s$. This implies, $\chi^{'}(G_{0,s}) = s$ and thereby, $\chi^{'}(G_{0,s}) = s = \chi^{'}_{\mathbb{Z}-sum}(G_{0,s})$.

\par
\noindent
\item [\rm {(c)}]    Let $n \geq 2$, $V(S_n)$ = $\{u_1, v_1, v_2, . . . , v_{n-1}\}$, $d(u_1) = n-1$ and $d(v_1)$ = 1 = $d(v_2)$ = . . . = $d(v_{n-1})$. $S_n$ is an integral sum graph \cite{ch90}; $S_n$ has $n-1$ edge-sum classes, each a singleton set and thereby $\chi^{'}_{\mathbb{Z}-sum}(S_n)$ = $n-1$ and $\triangle (S_n) = d(u_1) = n-1 = \chi^{'}(S_n)$ since all the $n-1$ edges at $u_1$ take different colors. 

\par
\noindent
\item [\rm {(d)}]    For $s \geq 2$, the edge-sum classes of $G_{-1,s}$ are $E_0, E_1,\ldots , E_s, E_{-1}$, each one is non-empty and hence $\chi^{'}_{\mathbb{Z}-sum}(G_{-1,s}) = s+2.$ Also, $\Delta (G_{-1,s}) = 1+s = deg(0) = deg(-1)$ and 

$\{ \{(0,1),$ $(-1,2)\},$ 

$\{(0,2),$ $(-1,3)\},$ 

$\{(0,3),$ ~ $(-1,4),$ $(1,2)\},$ 

$\{(0,4), (-1,5), (1,3)\},$~ 

$\{(0,5),$ $(-1,6),$ $(1, 4),$ $(2,3)\},$ 

$\{(0,6),$ $(-1, 7),$ $(1,5),$ $(2,4)\},$ 

$\{(0, 7),$ $(-1,8),$ $(1,6),$ $(2,5), (3, 4)\},$ 

$\{(0,8), (-1,9), (1,7), (2,6), (3, 5) \},$ 

. . . , 

$\{(0,s-1), (-1,s), (1,s-2),$ $(2,s-3),$ . . . , $(\left\lfloor \frac{s-3}{2}\right\rfloor, \left\lfloor \frac{s+2}{2}\right\rfloor )\},$ 

$\{(0, s),$ $(-1,1)\},$ 

$\{(0,-1),$ $(1,s-1),$ $(2,s-2),$ $(3,s-3),$ . . . , $( \left\lfloor \frac{s-1}{2}\right\rfloor,$ $\left\lfloor \frac{s+2}{2}\right\rfloor )\}\}$
\\
is a set of edge color classes of $G_{-1,s}$ which is of order $s+1$. This implies, edge color number of $G_{-1,s}$ is $s+1$, using Vizing's theorem \cite{ha69}. This implies, for $s \geq 2,$ $\chi^{'}(G_{-1,s}) = 1+s \neq \chi^{'}_{\mathbb{Z}-sum}(G_{-1,s}) = s+2$. 

 \par
\noindent
\item [\rm {(e)}]  At first, a set of edge color class of order $2s$ of $G_{-s,s}$ is obtained for each $s$, $s$ = 2,3,4,5,6 as follows:
	\begin{enumerate}
			\item [\rm (i)] A set of edge color classes of order 4 of $G_{-2, 2}$ is $\{\{(0, 2),$ $(-1, 1)\},$ $\{(0, 1),~ (-2,~2)\},$ $\{(0,~ -1), ~(-2, ~1)\},$ $\{(0,~-2),$ $(-1,2)\}\}$. 
				\item [\rm (ii)] A set of edge color classes of order 6 of $G_{-3,3}$ is $\{\{(0,3),$ $(1,$ $2)\}, \{(0,2), (-3,1)\},$ $\{(0,1),(-2,3), (-1,2)\}, \{(0,-1),$ $(-3,$ $2), (-2,1)\}, \{(0,-2), (-1,3)\},$ $\{(0,-3),$ $(-2,-1)\}\}$. 
				\item [\rm (iii)] A set of edge color classes of order 8 of $G_{-4,~4}$ is $\{\{(0,~4),$ $(-4,3), (-2,2), (-1,1)\}$, $\{(0,~3), (-4,~4), (1,~2)\},$ $\{(0,~2),$ $(-4,1), (-3,4), (-1,3)\},$ $\{(0,1), (-3,~3), (-2,~4), (-1,~2)\},$ $\{(0,~-1), (-3,~2), (-2,~3)\},$ $\{(0,~-2), (-3,~1), (-1, 4)\},$ $\{(0,-3),$ $(-4,2),$ $(-2,-1), (1,~3)\},$ $\{(0,~-4), (-3,~-1),$ $(-2,1)\} \}$. 
				\item [\rm (iv)] A set of edge color classes of order 10 of $G_{-5,5}$ is $\{\{(0,~5),$ $(-5,1), (-4,-1), (-3,-2),(2,3)\},$ $\{(0,4), (-5,3), (-4,1),$ $(-3,2), (-2,5)\},$ $\{(0,3), (-1,~4), (1,~2)\},$ $\{(0,~2), (-3,~4),$ $(-2,1), (-1,3)\},$ $\{(0,1), (-4,4), (-3,5), (-2,~3), (-1,~2)\},$ 
$\{(0,~-1), (-5,~5), (-4,~3), (-3,~1), (-2,~2)\},$ $\{(0,~-2),$ $(-5, 4), (-4,~2), (-3,~3), (-1,~1)\},$ $\{(0,~-3), (-2,~-1),$ $(1,3)\},$ $\{(0, -4), (-5,2), (-2,4), (-1,5)\},$ $\{(0,-5),$ $(-4,5),$ $(-3,-1), (1, 4)\}\}$. 
				
\item [\rm (v)] A set of edge color classes of order 12 of $G_{-6,6}$ is $\{\{(0,6),$ $(-4, ~3), (2,~4), (-1,~1)\},$ $\{(0,~5), (-1,~6), (1,~4), (2,~3)\},$ $\{(0,4),$ $(-6, 1), (-4,2), (-2,~6), (-1,~ 5)\},$ $\{(0,~3),$  $(-6,~6),$ $(-4,5), (-3,1), (-2,2), (-1,4)\},$ $\{(0,~2), (-6,~4), (-5, ~1),$ $(-3, 3),$ $(-2,5) \},$ $\{(0,1), (-5,6), (-4,~4), (-3,~5), (-2,~3),$ $(-1, 2)\},$ $\{(0, ~-1), (-5,~5),$ $(-3,~2), (-2,~1) \},$ $\{ (0,~-2),$ $(-5,~4), (-4, ~1), ~(-3,~6), ~(-1,~3)\},$ $\{(0,~-3), (-6,~3),$ $(-5,~-1), (-4,~-2), (1,~5)\},$ $\{(0, ~-4), (-6,~5), (-5,~2),$ $(-3,~-1), (-2,~4), (1,~3)\},$ $\{(0,~-5), (-6,~2), (-4,~-1),$ $(-3,-2)\},$ $\{(0,-~6), (-5, ~3), (-4,~6), (-3,~ 5), (-2,-~1),$ $(1,2)\}\}$.
\end{enumerate} 
			
			Using $(a)$,~ $\chi^{'}_{\mathbb{Z}-sum}(G_{r,s}) = -r+s+1$ for $-r,s \in \mathbb{N}$ and thereby $\chi^{'}_{\mathbb{Z}-sum}(G_{-s,s}) = 2s+1$ for $s \in \mathbb{N}$. For $s \geq 2,$ $\Delta (G_{-s,s}) = 2s = deg(0)$. This implies, for $s = 2,3,4,5,6,$ the edge color number of $G_{-s,s}$ is $2s$ using Vizing's theorem \cite{ha69}. Thus, for $2 \leq s \leq 6,$ $\chi^{'}(G_{-s,s}) = 2s \neq \chi^{'}_{\mathbb{Z}-sum}(G_{-s,s}) = 2s+1$. 
\end{enumerate}
Hence the theorem is proved. 	
\end{proof}

In Theorem \ref{c5}, we could prove that $\chi^{'}(G_{-s,s})  \neq \chi^{'}_{\mathbb{Z}-sum}(G_{-s,s})$ for $2 \leq s \leq 6$ and also $\chi^{'}(G_{-1,t})  \neq \chi^{'}_{\mathbb{Z}-sum}(G_{-1,t})$ for $t \geq 2$. In the following theorem, we prove its general result on $G_{r, s}$.
		
\begin{theorem}\quad {\rm \label{c6} For $s \geq 2$, $\chi^{'}(G_{r, s})$ = $|r| + s$ and $\chi^{'}_{\mathbb{Z}-sum}(G_{r, s})$ = $|r| + s +1$, $r < 0 < s$ and $-r \leq s$. And in particular, $\chi^{'}(G_{-s, s})$ = $2s$ and $\chi^{'}_{\mathbb{Z}-sum}(G_{-s, s})$ = $2s+1$.} 
\end{theorem}
\begin{proof} \quad In Theorem \ref{c5}, we proved that 
$\chi^{'}_{\mathbb{Z}-sum}(G_{r,s}) = |r| + s + 1$ for $-r,s \in \mathbb{N}$. Moreover, $\triangle (G_{r, s})$ = $-r+s$ and using Vizig's Theorem \cite{ha69}, $\triangle (G_{r, s})$ $\leq$ ${\chi}^{'}(G_{r,s}) \leq \triangle (G_{r, s})+1$.  

\vspace{.2cm}
\noindent
{\bf Claim.}  ${\chi}^{'}(G_{r,s}) = -r + s$ for $s \geq 2$, $r < 0 < s$ and $-r \leq s$.  

Here, we assign a coloring to the edges of $G_{r, s}$ with $-r+s$ colors and prove that it is a proper edge coloring of $G_{r, s}$ and so ${\chi}^{'}(G_{r,s}) = -r + s$. 

Let $c_j$ denote $j^{th}$ color assigned to an edge and $C_k$ denote the color class of edges each with color $k$ in $G_{r,s}$, $r < 0 < s$ and $j,k,-r,s\in\mathbb{N}$.   

Color the edges of $G_{r, s}$ as follows for $s \geq 2$, $r < 0 < s$, $-r \leq s$ and $i,j,k,-r,s\in\mathbb{N}$.

\vspace{.2cm}
\hspace{1.1cm}  $(0, j)$ $\mapsto$ $c_j$, ~$1 \leq j \leq s$;

\vspace{.2cm}
\hspace{.85cm}  
$(0, -i)$ $\mapsto$ $c_{i+s}$, ~$1 \leq i \leq |r|$;

\vspace{.2cm}
\hspace{.85cm}  
$(-i, j)$ $\mapsto$ $c_{i+j}$, ~ $1 \leq i \leq |r|$ and $1 \leq j \leq s-1$;

\vspace{.2cm}
\hspace{.85cm}  
$(-i, s)$ $\mapsto$ $c_{2i+s}$,~ $1 \leq i \leq \left\lfloor \frac{|r|}{2}\right\rfloor$;

\vspace{.2cm}
\hspace{.85cm}  
$(-i, s)$ $\mapsto$ $c_{i-\left\lfloor \frac{|r|}{2}\right\rfloor}$, ~$\left\lfloor \frac{|r|}{2}\right\rfloor < i \leq |r|$;

\vspace{.2cm}
\hspace{1.2cm}  
$(i, j)$ $\mapsto$ $c_{i+j+|r|}$,~ $1 \leq i,j,i+j \leq s$ and $i < j$ and

\vspace{.2cm}
\hspace{.5cm}  
$(-i, -j)$ $\mapsto$ $c_{i+j+s}$,~ $1 \leq i,j,i+j \leq |r|$ and $i < j$.

\vspace{.2cm}

It is clear from the above edge coloring that colors $c_1$ to $c_{s+|r|}$ are assigned to the edges of $G_{r,s}$ and no more edge colors. And also colors of edges at each vertex of $G_{r,s}$ are all distinct by the following. 

By Theorem \ref{a1}, we have $G_{r,s}$ $\cong$ $K_1*G_{-r}*G_s$. Also, $G_{-r}*G_s$ $\cong$ $G_{-r} \cup G_s \cup K_{-r, s}$ and vertex sets of $K_1$, $K_{r, s}$, $G_s$, $G_{-r}$ in $G_{r, s}$ are $\{0\}$, $\{1,2,...,s, -1,-2,...,r\}$, $\{1,2,...,s\}$, $\{-1,-2,...,r\}$, respectively. We show that the above edge coloring is a proper coloring by presenting the verification with respect to $K_1$, $K_{r, s}$, $G_s$, $G_{-r}$ in $G_{r, s}$ one by one.

Edges (incident) at 0 of $K_1$ in $G_{r,s}$ take all the $s+|r|$ colors $c_1$, $c_2$, $. . . ,$ $c_{s+|r|}$ and so the edge coloring is proper at 0 of $K_1$.

In $K_{-r,s}$, 
\begin{enumerate}
\item for $1 \leq j \leq s-1$, at $j$, colors of the $|r|$ edges are $c_{j+1}$, $c_{j+2}$, . . . , $c_{j+|r|}$ (which are distinct and also different from color $c_j$ of the edge $(0, j)$ at $j$) and thus the edge colors in this case are $c_2, c_3, . . . , c_{s-1+|r|}$. 

\item at $s$, colors of the $|r|$ edges are $c_{2\left\lfloor \frac{|r|}{2}\right\rfloor+s}$, $c_{2(\left\lfloor \frac{|r|}{2}\right\rfloor-1)+s}$, ..., $c_{2+s}$, $c_1, c_2, . . . , c_{|r|-\left\lfloor \frac{|r|}{2}\right\rfloor}$ which are distinct and also different from color $c_s$ of the edge $(0, s)$ and $2\left\lfloor \frac{|r|}{2}\right\rfloor+s$ $\leq$ $|r|+s$. 

\item for $1 \leq i \leq \left\lfloor \frac{|r|}{2}\right\rfloor$, at $-i$, colors of the $s$ edges $(-i, j)$ are $c_{i+1}$, $c_{i+2}$, . . . , $c_{i+s-1}$, $c_{2i+s}$ which are distinct and also different from color $c_{i+s}$ of the edge $(0, -i)$, $1 \leq j \leq s$.  

\item for $\left\lfloor \frac{|r|}{2}\right\rfloor < i \leq |r|$, at $-i$, colors of the $s$ edges $(-i, j)$ are $c_{i+1}$, $c_{i+2}$, . . . , $c_{i+s-1}$, $c_{i-\left\lfloor \frac{|r|}{2}\right\rfloor}$ which are distinct and also  different from color $c_{i+s}$ of the edge $(0, -i)$, $1 \leq j \leq s$.

Thus the edge coloring is proper in this case.
\end{enumerate}

In $G_{s}$, for $1 \leq j \leq s$, colors of edges are $c_{j+i+|r|}$ where $1 \leq i \leq s$, $i \neq j$ and $3 \leq i+j \leq s$ and so $c_{3+|r|}$, $c_{4+|r|}$, . . . , $c_{s+|r|}$ are the edge colors in this case. Edge colors at $i$ are $c_{j+i+|r|}$ which are distinct and also different from $c_i$ of the edge $(0, i)$ and also of other possible colors of edges which incident at $i$ in $G_{r,s}$ where $j \neq i$ and $3 \leq i+j \leq s$. The same arguments also hold when $i$ and $j$ are interchanged.

In $G_{-r}$, for $1 \leq i \leq |r|$, colors of edges are $c_{i+j+s}$ where $1 \leq j \leq |r|$, $i \neq j$ and $3 \leq i+j \leq |r|$ and so $c_{3+s}$, $c_{4+s}$, . . . , $c_{|r|+s}$ are the edge colors in this case. Edge colors at $-i$ are $c_{j+i+s}$ which are distinct and also different from $c_{i+s}$ of the edge $(0, -i)$ and also of other possible colors of edges which incident at $-i$ in $G_{r,s}$ where $j \neq i$ and $3 \leq i+j \leq |r|$. The same arguments also hold when $-i$ and $-j$ are interchanged.

Thus, in the above edge coloring, $G_{r,s}$ takes only $|r|+s$ number of colors and colors of edges incident at each vertex of $G_{r,s}$ are all distinct and thereby ${\chi}^{'}(G_{r,s})$ = $|r|+s$ since $\triangle (G_{r, s})$ = $|r|+s$ and by Vizig's theorem \cite{ha69}, $\triangle (G_{r, s})$ $\leq$ ${\chi}^{'}(G_{r,s}) \leq \triangle (G_{r, s})+1$. Hence we get the result.

Figures 6 and 7 show the edge coloring, as given in this theorem, of the integral sum graphs $G_{-4, 7}$ and $G_{-5,7}$. In Figure 7 (as well as in Figure 6), $c_1$, $c_2$, . . . , $c_{12}$ take the colors gray, orange, cyan, purple, blue, green, brown, violet, magenta, red, olive, black, respectively. See Fig. 6 and Fig. 7. 
\end{proof}

	\begin{center}
\begin{figure}
	\begin{tikzpicture}
		\node (c0) at (9.5,5)  [circle,draw,scale=0.8]{0};
		\node (c1) at (13,4.75)  [circle,draw,scale=0.8]{1};
		\node (c2) at (15.75,3.25)  [circle,draw,scale=0.8] {2};
		\node (c3) at (17,1.5)  [circle,draw,scale=0.8] {3};
		\node (c4) at (16.5,-.5)  [circle,draw,scale=0.8] {4};
		\node (c5) at (15,-2)  [circle,draw,scale=0.8]{5};
		\node (c6) at (13,-3) [circle,draw,scale=0.8] {6};
		\node (c7) at (11,-3.75) [circle,draw,scale=0.8] {7};
		\node (c8) at (7,3.25) [circle,draw,scale=0.8] {-1};
		\node (c9) at (6,.75) [circle,draw,scale=0.8] {-2};
		\node (c10) at (6.5,-1)  [circle,draw,scale=0.8] {-3};
		\node (c11) at (7.5,-2.25)  [circle,draw,scale=0.8] {-4};
		
	\draw (c0)[gray, thick] --node [above] {$c_1$} (c1);
	\draw (c0)[orange, thick] --node[near start] {$c_2$}  (c2);
	\draw (c0)[cyan, thick] --node[near start] {$c_3$}  (c3);
	\draw (c0)[purple, thick] --node[ near start] {$c_4$}  (c4);
	\draw (c0)[blue, thick] --node[ near start] {$c_5$} (c5);
	\draw (c0)[green, thick]--node[near start] {$c_6$}  (c6);
	\draw (c0)[brown, thick] --node[near start] {$c_7$}  (c7);
	\draw (c0)[violet, thick] --node[near start] {$c_{8}$}  (c8);
	\draw (c0)[magenta, thick] --node[near start] {$c_{9}$} (c9);
	\draw (c0)[red, thick] --node[near start] {$c_{10}$}  (c10);
	\draw (c0)[olive, thick] --node[near start] {$c_{11}$}  (c11);
		
	\draw (c1)[orange, thick] --node[near start] {$c_2$} (c8);
	\draw (c1)[cyan, thick] --node[near start] {$c_3$}(c9);
	\draw (c1)[purple, thick] --node[near start] {$c_4$} (c10);
	\draw (c1)[blue, thick] --node[near start] {$c_5$}  (c11);
		
	\draw (c2)[cyan, thick] --node[near start] {$c_3$} (c8);
  \draw (c2)[purple, thick] --node[near start, right] {$c_4$} (c9);
  \draw (c2)[blue, thick] --node[near start] {$c_5$}(c10);
	\draw (c2)[green, thick] --node[near start] {$c_6$} (c11);
	
   \draw (c3)[purple, thick] --node[near start] {$c_4$} (c8);
	 \draw (c3)[blue, thick] --node[near start] {$c_5$}(c9);
	 \draw (c3)[green, thick] --node[near start] {$c_6$} (c10);
	 \draw (c3)[brown, thick] --node[near start] {$c_7$} (c11);
	
	\draw (c4)[blue, thick] --node[near start] {$c_5$}(c8);
	\draw (c4)[green, thick] --node[near start] {$c_6$}(c9);
	\draw (c4)[brown, thick] --node[near start] {$c_7$}  (c10);
	\draw (c4)[violet, thick] --node[near start] {$c_{8}$} (c11);

	\draw (c5)[green, thick] --node[near start] {$c_6$} (c8);
	\draw (c5)[brown, thick] --node[near start] {$c_7$} (c9);
	\draw (c5)[violet, thick] --node[near start] {$c_{8}$}(c10);
  \draw (c5)[magenta, thick] --node[near start] {$c_{9}$} (c11);

	\draw (c6)[brown, thick] --node[near start] {$c_7$} (c8);
	\draw (c6)[violet, thick] --node[near start] {$c_{8}$} (c9);
	\draw (c6)[magenta, thick] --node[near start] {$c_{9}$} (c10);
	\draw (c6)[red, thick] --node[near start] {$c_{10}$} (c11);
		
	\draw (c7)[magenta, thick] --node[near start] {$c_{9}$} (c8);
	\draw (c7)[olive, thick] --node[near start] {$c_{11}$} (c9);
	
	\draw (c7)[gray, thick] --node[near start] {$c_{1}$} (c10);
	\draw (c7)[orange, thick] --node[below] {$c_2$} (c11);
	
	\draw (c1)[brown, thick] --node[above] {$c_7$}  (c2);
	\draw (c1)[violet, thick] --node[near start] {$c_{8}$}(c3);
	\draw (c1)[magenta, thick] --node[near start] {$c_{9}$} (c4);
	\draw (c1)[red, thick] --node[near start] {$c_{10}$} (c5);
	\draw (c1)[olive, thick] --node[near start] {$c_{11}$}(c6);
	
	\draw (c2)[magenta, thick] --node[right] {$c_{9}$} (c3);
	\draw (c2)[red, thick] --node[near start] {$c_{10}$} (c4);
	\draw (c2)[olive, thick] --node[near start] {$c_{11}$}(c5);
	
	\draw (c3)[olive, thick] --node[right] {$c_{11}$}(c4);
	
	\draw (c8)[red, thick] --node[left] {$c_{10}$} (c9);
  \draw (c8)[olive, thick] --node[near start] {$c_{11}$}(c10);
	\end{tikzpicture}
Fig. 6. $G_{-4,7}$ with edge coloring.
\end{figure}    
\end{center}
\begin{center}
	\begin{tikzpicture}
		\node (c0) at (10,6.5)  [circle,draw,scale=0.8]{0};
		\node (c1) at (14.5,6)  [circle,draw,scale=0.8]{1};
		\node (c2) at (16.5,3.75)  [circle,draw,scale=0.8] {2};
		\node (c3) at (17,1.5)  [circle,draw,scale=0.8] {3};
		\node (c4) at (16,-.5)  [circle,draw,scale=0.8] {4};
		\node (c5) at (14.5,-2)  [circle,draw,scale=0.8]{5};
		\node (c6) at (13,-3) [circle,draw,scale=0.8] {6};
		\node (c7) at (11,-3.5) [circle,draw,scale=0.8] {7};
		
		\node (c8) at (7,5) [circle,draw,scale=0.8] {-1};
		\node (c9) at (6,3) [circle,draw,scale=0.8] {-2};
		\node (c10) at (6,1)  [circle,draw,scale=0.8] {-3};
		\node (c11) at (6.5,-.5)  [circle,draw,scale=0.8] {-4};
		\node (c12) at (7.25,-2) [circle,draw,scale=0.8] {-5};
		
	\draw (c0)[gray, thick] --node[][above]{$c_1$} (c1);
	\draw (c0)[orange, thick] --node[near start] {$c_2$}  (c2);
	\draw (c0)[cyan, thick] --node[near start] {$c_3$}  (c3);
	\draw (c0)[purple, thick] --node[near start] {$c_4$}  (c4);
	\draw (c0)[blue, thick] --node[near start] {$c_5$}  (c5);
	\draw (c0)[green, thick]--node[near start] {$c_6$}  (c6);
	\draw (c0)[brown, thick] --node[near start] {$c_7$}  (c7);
	\draw (c0)[violet, thick] --node[near start] {$c_{8}$}  (c8);
	\draw (c0)[magenta, thick] --node[near start] {$c_{9}$}  (c9);
	\draw (c0)[red, thick] --node[near start] {$c_{10}$}  (c10);
	\draw (c0)[olive, thick] --node[near start] {$c_{11}$}  (c11);
\draw (c0)[black, thick] --node[near start] {$c_{12}$} (c12);
			
		\draw (c1)[orange, thick] --node[near start] {$c_2$} (c8);
		\draw (c1)[cyan, thick] --node[near start] {$c_3$}  (c9);
		\draw (c1)[purple, thick] --node[near start] {$c_4$} (c10);
		\draw (c1)[blue, thick] --node[near start] {$c_{5}$}(c11);
		\draw (c1)[green, thick] --node[near start] {$c_{6}$} (c12);
		
		\draw (c2)[cyan, thick] --node[near start] {$c_3$}  (c8);
		\draw (c2)[purple, thick] --node[near start] {$c_4$} (c9);
		\draw (c2)[blue, thick] --node[near start] {$c_{5}$}(c10);
		\draw (c2)[green, thick] --node[near start] {$c_{6}$} (c11);
		\draw (c2)[brown, thick] --node[near start] {$c_{7}$} (c12);
		
		\draw (c3)[purple, thick] --node[near start] {$c_4$} (c8);
		\draw (c3)[blue, thick] --node[near start] {$c_{5}$} (c9);
		\draw (c3)[green, thick] --node[near start] {$c_{6}$} (c10);
		\draw (c3)[brown, thick] --node[near start] {$c_{7}$}(c11);
		\draw (c3)[violet, thick] --node[near start] {$c_{8}$} (c12);
		
		\draw (c4)[blue, thick] --node[near start] {$c_{5}$} (c8);
		\draw (c4)[green, thick] --node[near start] {$c_{6}$} (c9);
		\draw (c4)[brown, thick] --node[near start] {$c_7$}  (c10);
		\draw (c4)[violet, thick] --node[near start] {$c_8$}(c11);
		\draw (c4)[magenta, thick] --node[near start] {$c_9$}(c12);
		
		\draw (c5)[green, thick] --node[near start] {$c_{6}$} (c8);
		\draw (c5)[brown, thick] --node[near start] {$c_{7}$} (c9);
		\draw (c5)[violet, thick] --node[near start] {$c_{8}$}(c10);
		\draw (c5)[magenta, thick] --node[near start] {$c_9$} (c11);
		\draw (c5)[red, thick] --node[near start] {$c_{10}$} (c12);
		
		\draw (c6)[brown, thick] --node[near start] {$c_{7}$} (c8);
		\draw (c6)[violet, thick] --node[near start] {$c_{8}$} (c9);
	\draw (c6)[magenta, thick] --node[near start] {$c_{9}$} (c10);
		\draw (c6)[red, thick] --node[near start] {$c_{10}$} (c11);
	\draw (c6)[olive, thick] --node[near start] {$c_{11}$} (c12);
		
		\draw(c7)[magenta, thick] --node[near start] {$c_9$} (c8);
		\draw (c7)[olive, thick] --node[near start] {$c_{11}$} (c9);
		
		\draw (c7)[gray, thick] --node[near start] {$c_1$} (c10);
	\draw (c7)[orange, thick] --node[near start] {$c_{2}$} (c11);
\draw (c7)[cyan, thick] --node[near start][left] {$c_{3}$} (c12);

  	\draw (c1)[violet, thick] --node[above] {$c_8$}  (c2);
		\draw (c1)[magenta, thick] --node[near start] {$c_9$}(c3);
		\draw (c1)[red, thick] --node[near start] {$c_{10}$} (c4);
		\draw (c1)[olive, thick] --node[near start] {$c_{11}$} (c5);
	  \draw (c1)[black, thick] --node[near start] {$c_{12}$}(c6);
	
		\draw (c2)[red, thick] --node[right] {$c_{10}$} (c3);
		\draw (c2)[olive, thick] --node[near end] {$c_{11}$} (c4);
	  \draw (c2)[black, thick] --node[near start] {$c_{12}$}(c5);
		
		\draw (c3)[black, thick] --node[right] {$c_{12}$}(c4);
		
		\draw (c8)[red, thick] --node[left] {$c_{10}$}(c9);
		\draw (c8)[olive, thick] --node[near end] {$c_{11}$}(c10);
   \draw (c8)[black, thick] --node[near start] {$c_{12}$} (c11);

	\draw (c9)[olive, thick] --node[left] {$c_{12}$} (c10);
	\end{tikzpicture}	
	Fig. 7. $G_{-5, 7}$ with edge coloring.
\end{center}

\noindent 
\textbf{Open Problem.}		Characterize those integral sum graphs $G^+(S)$ for which
\begin{enumerate}
\item  	$\chi^{'}_{\mathbb{Z}-sum}(G^+(S))$ = $|G^+(S)|$. \\
\item   $\chi^{'}(G^+(S))$ = $\chi^{'}_{\mathbb{Z}-sum}(G^+(S))$.\\
\item   $\chi^{'}(G^+(S)) < \chi^{'}_{\mathbb{Z}-sum}(G^+(S))$. \\
\item  $\chi^{'}(G^+(S)) > \chi^{'}_{\mathbb{Z}-sum}(G^+(S)).$ \\
\item $G^+(S)$ is neither of (2) nor of (3) and (4). (Check whether such graphs exist.)
\end{enumerate}

\vspace{.5cm}
\noindent
\textbf{Conflict of Interest}

\noindent
\textit {The authors declare that there is no conflict of interests regarding the publication of this paper.}

\vspace{.5cm}
\noindent
\textbf{Acknowledgement}

\noindent
\textit {The first author express his sincere thanks to the Central University of Kerala, Kasaragod - 671 316, Kerala, India and St.Jude's College, Thoothoor - 629 176, Tamil Nadu, India for providing facilities to carry out this research work.}

\end{document}